\documentclass[12pt,a4paper]{article}
\usepackage{graphicx,txfonts}

\begin{document}

\title{\bf Obituary: \,  Walter Benz \,  (1931--2017)}
\author{Alexander Kreuzer \and Hans Havlicek}
\date{}
\maketitle

\begin{figure}[h!]
  \centering
  \includegraphics[scale=0.3]{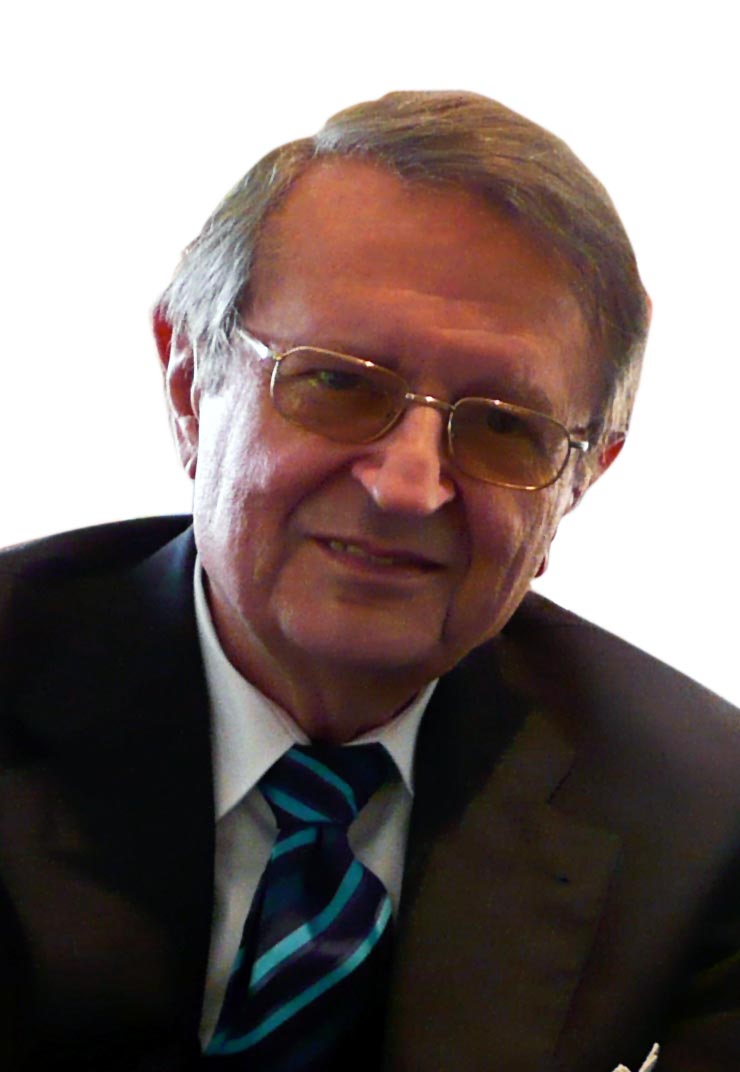}
	\bigskip
	\center{Walter Benz}
	\label{}
\end{figure}

Professor Dr.\ Dr.\ h.c.\ Walter Benz, an internationally renowned, highly
esteemed and influential mathematician, passed away on 13 January 2017. His
research enriched geometry and he founded this journal with Birkh\"auser in
1971. We have lost a very active researcher, a kind colleague and a reliable
friend. We will always keep him in fond memory and express our deepest sympathy
to his family and all who knew him.

\par

His live can be described accurately with the following statement, which he
made in 2007 in an interview given to Uta Hartmann (cf. [1]):

\par\smallskip

\begin{center} \parbox{9.8 cm}{``Die Mathematik war mein Leben,
ist mein Leben. Meine Begeisterung f\"ur dieses Fach ist ungebrochen. Es gibt
keinen Tag an dem ich mich nicht auch ernsthaft mit Mathe\-matik
besch\"aftige.''}
\end{center}

\par\smallskip

Walter Benz was born in Lahnstein, Germany, on 2 May 1931. During the last days
of the second world war his uncle taught and familiarized him with the books
``Differentialgleichungen'' by Karl-Heinrich Weise and
``Differentialgeometrie'' by Wolfgang Haak. These books founded the lifelong
bond of love of Walter Benz to mathematics. He started his study of mathematics
and physics at the Johannes-Gutenberg-Universit\"at Mainz in 1951, completed
his ``Staatsexamen'' in 1955 and received his Ph.D. one year later under the
supervision of Robert Furch (1894--1967) with a thesis in geometry, entitled
``Axiomatischer Aufbau der Kreisgeometrie auf Grund von
Doppelverh\"altnissen''. In the year 1959 Walter Benz qualified as a university
lecturer with the thesis ``\"Uber Winkel- und Transitivit\"atseigenschaften in
Kreisebenen''.

\par

In the following years Walter Benz was teaching mathematics in Frankfurt
(1961--1965), Mainz (1961), W\"urzburg (1963--64), Bochum (1966--1974) and at
the University of Waterloo (1966--67, 1969), since 1966 as a full professor,
and from 1971 to 1974 as the Dean of the Mathematical Faculty of the University of
Bochum. Then, in 1974, he was appointed as full professor at the University of
Hamburg, where he also served as the Dean of the Mathematical Faculty
(1985--1988), and where he stayed until his retirement in 1997. He also was
visiting professor at Bologna, Brescia, Freiburg, Kuwait, Rome, Thessaloniki
and Trieste.

\par

Walter Benz is well-known for his comprehensive work in several parts of
geometry. He was researching in the geometries of M\"obius, Laguerre, Minkowski
and Lie, contributed deep results to the discipline ``characterizations of
geometrical mappings under mild hypotheses'' and wrote papers about functional
equations. Furthermore, he was involved with hyperbolic geometry, metric
geometry and special relativity. He published more than 180 articles and
several books, for example the highly influential early book ``Vorlesungen
\"uber Geometrie der Algebren''. The planes of M\"obius, Laguerre and Minkowski
are known today under the name ``Benz planes''.

\par

Walter Benz's achievements have been recognized by an honorary doctorate from
the University of Sofia (1995), an honorary fellowship in the ``Mathematische
Gesellschaft in Hamburg'' (1996), and the honorary fellowship (2004) of the
Hungarian Academy of Sciences. Furthermore he was a member of the European
Academy of Sciences (2004). He was an editor of the journals ``Jahresberichte
der DMV'', ``Journal of Geometry'', ``Aequationes Mathematicae'',
``Abhandlungen aus dem Mathe\-matischen Seminar der Universit\"at Hamburg'',
``Results in Mathematics'', ``Atti del Seminario Matematico e Fisico
dell'Universit\`a di Modena'' and ``Mitteilungen der Mathematischen
Gesellschaft in Hamburg''.

\par

His inexhaustible knowledge about other famous mathematicians was very
impressing. He knew not only their teachers and pupils, but also was able to
tell humorous and distinguishing anecdotes.

\par

Walter Benz will stay alive by his mathematical work and will be always
remembered by his friends and colleagues.

\bibliographystyle{plain}

\small\noindent
Alexander Kreuzer\\
Fachbereich Mathematik\\
Universit\"at Hamburg\\
Bundesstra{\ss}e 55\\
D-20146 Hamburg\\
Germany\\
\texttt{kreuzer@math.uni-hamburg.de}

\par\medskip\noindent
Hans Havlicek\\
Institut f\"ur Diskrete Mathematik und Geometrie\\
Technische Universit\"at
Wien \\ Wiedner Hauptstra{\ss}e 8--10\\
A-1040 Wien\\
Austria\\
\texttt{havlicek@geometrie.tuwien.ac.at}

\end{document}